\input amstex
\documentstyle{amsppt}
\input bull-ppt
\keyedby{rio1/amh}
\topmatter
\cvol{31}
\cvolyear{1994}
\cmonth{October}
\cyear{1994}
\cvolno{2}
\cpgs{208-212}
\title Singular continuous spectrum is generic\endtitle
\author R. del Rio, S. Jitomirskaya, N. 
Makarov, and B. Simon\endauthor
\shortauthor{del Rio, Jitomirskaya, Makarov, and Simon}
\address {\rm (R. del Rio, N. Makarov, and B. Simon)}
Division of Physics, Mathematics and Astronomy, 253-37, 
California Institute of Technology, Pasadena, California  
91125\endaddress 
\cu\nofrills Current address, R.del Rio: 
IIMAS-UNAM, Apdo.~Postal 20-726, 
Admon No.~20, 01000 Mexico D.F., Mexico\endcu 
\ml\nofrills{\it E-mail address, \rm R. del Rio:}\enspace
\tt delrio\@redvax1.dgsca.unam.mx
\mlpar{\it E-mail address, \rm N. Makarov:}\enspace\tt 
makarov\@caltech.edu
\mlpar{\it E-mail address, \rm B. Simon:}\enspace\tt 
bsimon\@caltech.edu
\endml
\address {\rm (S. Jitomirskaya)} 
Department of Mathematics, University of California, 
Irvine, California 92717\endaddress
\ml szhitomi\@math.uci.edu \endml
\thanks The first author's research was partially 
supported by 
DGAPA-UNAM and CONACYT\endthanks 
\thanks This material is based upon work of the third 
author supported by 
National Science Foundation grant DMS-9207071. The 
government 
has certain rights in this material\endthanks
\thanks This material is based upon work of the first and 
fourth authors
 supported by  
National Science Foundation  grant DMS-9101715. The 
government 
has certain rights in this material\endthanks 
\date May 3, 1993\enddate
\abstract In a variety of contexts, we prove that singular 
continuous 
spectrum is generic in the sense that for certain natural 
complete metric 
spaces of operators, those with singular spectrum are a 
dense $G_\delta$.\endabstract
\endtopmatter

\document

In the spectral analysis of various operators of 
mathematical physics, 
a key step, often the hardest, is to prove that the 
operator has no 
continuous singular spectrum, that is, that the spectral 
measures for 
the operators have only pure point and absolutely 
continuous parts.  
Examples are the absence of such spectrum for $N$-body 
Schr\"odinger 
operators [3, 19] and for the one-dimensional random models 
[12, 7, 8, 24, 18].

Our goal here is to announce results that show that a 
singular 
continuous spectrum  lies quite close to many operators by 
proving 
it is often generic in the Baire sense.  Detailed proofs 
and further 
results will appear in three papers: one for general 
operators [22], 
one for rank-one perturbations [6], and one for 
almost-periodic 
Schr\"odinger operators [23].

Precursors of our results include work on generic ergodic 
processes 
[15, 21] and on special energies for Sch\"rodinger 
operators/Jacobi matrices 
[11, 4, 5].  Gordon [13, 14] independently (and 
presumably, before 
us) proved Theorem 5.  His method of proof is very 
different from ours.

Recall that the Baire category theorem implies that if $X$ 
is a 
complete metric space, a countable intersection of dense 
$G_\delta$ 
is still a dense $G_\delta$ and that if $X$ is perfect, 
then any dense 
$G_\delta$ has uncountable intersection with any open ball.

Our first two results are for one-body Schr\"odinger 
operators and for 
the ``generic Anderson model''.

\proclaim{Theorem 1 \rm [22]}  Let $C_{\infty}(\Bbb 
R^{\nu})$ denote the 
continuous functions on $\Bbb R^{\nu}$ vanishing at 
infinity in the $\|\cdot 
\|_{\infty}$ norm.  Then for a dense $G_\delta$ of $V\in 
C_{\infty}
(\Bbb R^{\nu})$, $-\Delta+V$ has purely singular 
continuous spectrum on 
$(0, \infty)$.
\endproclaim
\rem{ 
Remarks}  1.  If $V(x)=0(|x|^{-1-\epsilon})$ at infinity, 
it is 
known [1, 20] that $-\Delta+V$ has absolutely continuous 
spectrum on 
$(0, \infty)$ with a possible set of eigenvalues.

2.  There is a similar result [22] for $\{V\mid (1+
x^{2})^{\alpha/2} 
V\in C(\Bbb R^{\nu})\}$ with norm $|||V|||=\|(1+
x^{2})^{\alpha/2} 
V\|_{\infty}$ so long as $\alpha <\frac{1}{2}$.
\eject

3.  In one dimension, there is a similar result for Jacobi 
matrices [22].
\endrem 
\proclaim{Theorem 2 \rm[22]}  Let $a<b$ be fixed and let 
$\Omega=
[a, b]^{\Bbb Z^{\nu}}$, functions $v:\Bbb Z^{\nu}\to 
[a,b]$ with the 
\RM(compact metrizable\/\RM) Tychonoff topology.  Given 
any $v$, let $h(v)$ be the Jacobi 
matrix on $l^2 (\Bbb Z^{\nu})$ by 
$$
(h(v)u)(n)=\sum_{|j|=1}u(n+j)+v(n)u(n). \tag 1
$$
Then for a dense $G_\delta$ in $\Omega$, $h(v)$ has 
spectrum $[a-2\nu, 
b+2\nu]$, and the spectrum is purely singular continuous.
\endproclaim
\rem{ 
Remark}  If $\nu=1$ or if $\nu$ is arbitrary and $b-a$ is 
large, 
it is known that if $\Omega$ is given the product measure 
$\operatornamewithlimits{X}\limits_{j\in\Bbb Z^{\nu}} 
(b-a)^{-1}\, dx_{j}$, 
then for a.e.~$v$, $h(v)$ has only pure point spectrum [7, 
9, 8, 24, 2].  
So the generic Baire and generic Lebesgue behaviors are 
very different.
\endrem 

These results are not limited to Schr\"odinger operators.
\proclaim{Theorem 3 \rm[22]}  Let $X_a$ be the family of 
all self-adjoint 
operators, $A$, on a fixed separable Hilbert space, $\Cal 
H$, with 
$\|A\|\leq a$.  Give $X_a$ the metrizable topology of 
strong 
convergence.  Then $X_a$ has a dense $G_\delta$ of 
operators with 
$\operatorname{spec}(A)=[-a, a]$, and the spectrum is 
purely singular continuous.
\endproclaim
\proclaim{Theorem 4 \rm[22]}  Let $A$ be a fixed 
self-adjoint operator.  
Let $\Cal I_2$ be the Hilbert-Schmidt operators with 
Hilbert-Schmidt 
norm.  Then for a dense $G_\delta$ of $C$'s in $\Cal I_2$, 
the set of 
vectors
$$\align
&\{\psi\mid d\mu^{\psi}_{(A+C)}\text{ is purely singular 
continuous}\} \\
&\qquad \cup\{\psi\mid (A+C)\psi=E\psi;\ E\in
\operatorname{spec}_{\text{disc}}(A+C)\}
\endalign
$$
span $\Cal H$.
\endproclaim   
\rem{ 
Remarks}  1.  $d\mu_{D}^\psi$ is the spectral measure 
$(\psi, e^{isD}\psi)=\int e^{isE}\,d\mu_{D}^{\psi}(E)$.

2.  If $\text{spec}(A)$ is thin, for example, $A\equiv 0$, 
the vectors are 
the discrete eigenvectors.  But if $\text{spec}(A)$ 
contains an interval, 
$\text{spec}(A+C)$ will have lots of singular continuous 
spectrum.

3.  This is to be distinguished from the Weyl-von Neumann 
theorem 
[26, 27, 17] that there are $C$'s with $\|C\|_{2}$ 
arbitrarily small so that 
$A+C$ has only point spectrum.  Here we see that 
generically there 
will be singular continuous spectrum (for $A$ suitable).

4.  The result holds if $\Cal I_2$ is replaced by $\Cal 
I_p$ with 
$p>1$.  If $A$ has no a.c.~spectrum, it even holds for 
$\Cal I_1$.
\endrem 

These four theorems are rather soft with no hard 
estimates.  More 
subtle is the case of rank-one perturbations.  We will 
consider two 
closely related cases:
\roster
\item"{(a)}"  $A$ is a self-adjoint operator with cyclic 
vector $\varphi$; let 
$P$ be the projection onto $\varphi$, and let $A_\lambda 
=A+\lambda P$.
\item"{(b)}"  Let $H$ be the differential operator 
$-\frac{d^2}{dx^2}+V(x)$ on 
$[0, \infty)$ assumed to be limit point at infinity.  
$H_\theta$ is the 
self-adjoint operator with boundary condition $\cos\theta 
u(0)+\sin\theta 
u'(0)=0$.
\endroster
\proclaim{Theorem 5 \rm[6]}  {\rm (a)}  Suppose $A$ has an 
interval $[a, b]$ 
in its spectrum and the spectrum there has no 
a.c.~component.  Then
\roster
\item"\rom{(i)}"  There is a dense $G_\delta$, $C$, in 
$[a, b]$ so that if
$E\in C$, then $E$ is not an eigenvalue of any $A_\lambda$.
\item"\rom{(ii)}"  For a dense $G_\delta$, $L$, of $\Bbb 
R$, $A_\lambda$ has 
purely singular spectrum in $[a, b]$ if $\lambda\in L$.
\endroster

{\rm (b)}  Suppose for some $\theta _0$, $H_{\theta_0}$ 
has an interval 
$[a, b]$ in its spectrum and the spectrum there has no 
a.c.~component, then
\roster
\item"\rom{(i)}"  There is a dense $G_\delta$, $C$, in 
$[a, b]$ so that if
$E\in C$, then $E$ is not an eigenvalue of any $H_\theta$.
\item"\rom{(ii)}"  For a dense $G_\delta$, $L$, of $[0, 
\pi)$, 
$H_\theta$ has purely singular spectrum in $[a, b]$ if 
$\theta\in L$.
\endroster
\endproclaim
\rem{ 
Remarks}  1.  Case (i) implies that under the hypothesis 
of $E\in C$, 
either \break $\lim\limits _{x\to\infty}\frac{1}{x} %
\ln\|T_{E}(x)\|$ fails 
to exist or is $0$ when $T(E)$ is the fundamental matrix 
for the problem.  
This means that for many cases where one can only prove 
Lyapunov behavior for 
a.e.~$E$, there really is a set where the Lyapunov 
behavior fails [11, 4, 5].

2.  There are also results for general $A$ without any 
hypothesis on 
$\text{spec}(A)$ or absolute continuous spectrum.

3.  These results imply that in the Anderson model in the 
localized 
regime, varying $V(0)$ a little can produce singular 
spectrum.  Indeed, there 
are disjoint, locally uncountable sets with purely pure 
point spectrum when 
$V(0)$ is in one set and pure singular continuous spectrum 
when $V(0)$ is 
in the other set!
\endrem 

Another subtle class is the almost-periodic Schr\"odinger 
operators.  
We will consider functions $V$ on $\Bbb R$ or $\Bbb Z$ 
that are even and almost 
periodic (typical examples are $V(n)=\lambda\cos(\pi\alpha 
n)$ in the 
$\Bbb Z$ case and $V(x)=\lambda\cos(\pi x)+\mu\cos(\pi 
\alpha x)$ in 
the $\Bbb R$ case with $\alpha$ irrational) and define
$$\align
&\Bbb R\text{ case}\qquad H_\omega =-\frac{d^2}{dx^2}+
V_{\omega}(x)  \\
&\Bbb Z\text{ case}\qquad ({H_\omega}u)(n) =u(n+1)+u(n-1)+
V_{\omega}(n)u(n)
\endalign
$$
where $\omega$ is a point in the hull, $\Omega$, of $V$ 
and $V_\omega$ the 
corresponding potential (in the typical cases above, 
$\Omega=S^1$ and 
$S^{1}\times S^1$ with 
$V_{\theta}(n)=\lambda\cos(\pi\alpha n+\theta)$ and 
$V_{\theta,\psi}(x)=\lambda\cos(\pi x+\theta)+
\mu\cos(\pi\alpha x+\psi)$).  
$\Omega$ is a compact metric space in the Bohr topology.
\proclaim{Theorem 6 \rm[23]}  Let $V$ be an even 
almost-periodic 
potential on $\Bbb R$ or $\Bbb Z$.  Then\,\RM:
\roster
\item"\rom{(a)}"  For a dense $G_\delta$ in the hull, 
$H_\omega$ has no point 
spectrum.
\item"\rom{(b)}"  If for some point $\omega_0$ in the hull 
$H_{\omega_{0}}$ 
has no a.c.~spectrum, then for a dense $G_\delta$ in the 
hull, $H_\omega$ 
has purely singular continuous spectrum.
\endroster
\endproclaim
\ex{ 
Example \rm[23]}  In the $\Bbb Z$ case, if $V=\lambda\cos 
(\pi\alpha n+\theta)$ with $\lambda\geq 2$ and $\alpha$ 
irrational, then 
it follows that $H_\theta$ has purely singular spectrum 
for a dense 
$G_\delta$ of $\theta$.  When $\lambda$ is large [25, 10, 
16], it is known 
that we have pure point spectrum only for a set of 
$\theta$ of full 
Lebesgue measure.  Once again, we have locally uncountable 
sets of 
parameters with point spectrum for one parameter set and 
singular continuous 
spectrum in the other.
\endex 

\enddocument